\documentclass[11pt]{amsart}
\usepackage{amsfonts}
\usepackage{hyperref}
\usepackage[utf8]{inputenc}
\usepackage[T1]{fontenc}
\usepackage[dvips]{graphicx}
\usepackage{amssymb}
\usepackage{amsmath}
\usepackage{dsfont}
\usepackage{color}
\usepackage{latexsym}
\usepackage{bbm}
\usepackage{color}
\usepackage{amsthm}
\usepackage{multicol}
\usepackage[top   = 2.75cm,
bottom = 2.50cm,
left   = 2.50cm,
right  = 2.00cm]{geometry}

\bibstyle{plain}
\theoremstyle{plain}
\newtheorem{theorem}{Theorem}
\newtheorem{lemma}[theorem]{Lemat}

\newtheorem{corollary}[theorem]{Corollary}

\theoremstyle{remark}

\newtheorem{remark}[theorem]{Uwaga}

\begin{document}

\title[Singular distributions of random variables]{Singular distributions of random variables with independent digits of representation in numeral system with natural base and redundant alphabet
}
\author[M.\ V.\ Pratsiovytyi and 
S.\ P.\ Ratushniak]{ M.\ V.\ Pratsiovytyi and
S.\ P.\ Ratushniak}

\newcommand{\eacr}{\newline\indent}

\address{M.V. Pratsiovytyi\eacr
Institute of Mathematics of NASU,\eacr
 Dragomanov Ukrainian State University, \eacr Kyiv,
Ukraine\acr ORCID 0000-0001-6130-9413}
\email{prats4444@gmail.com}

\address{S.P. Ratushniak\eacr
Institute of Mathematics of NASU,\eacr
 Dragomanov Ukrainian State University,\eacr Kyiv,
Ukraine\acr ORCID 0009-0005-2849-6233}
\email{ratush404@gmail.com}

\keywords{ Singular distribution of a random variable, absolutely continuous distribution, infinite Bernoulli convolution governed by a series, numeral system with base~$3$ and a redundant digit, set of incomplete sums of a series, $N$-self-similar set, fractal dimension.\\This work was supported by a grant from the Simons
Foundation (SFI-PD-Ukraine-00014586, M.P., S.R\\
Published in Matematychni Studii, Vol. 63 (2) (2025), pp. 199–209.}

\date{\today}

\newcommand{\acr}{\newline\indent}

\begin{abstract}
Given natural parameters $s$ and $r$, where $2\leq s\leq r$, we consider the distribution of a random variable
$\xi=\sum\limits_{k=1}^{\infty}s^{-k}\xi_k\equiv\Delta^{r_s}_{\xi_1\xi_2...\xi_k...},$
where $(\xi_k)$ is a sequence of independent random variables taking values in $\{0,1,...,r\}$ with probabilities
$p_0,p_1,...,p_r$, respectively, and all $ p_i<1$.

In the case $s=3=r$, necessary and sufficient conditions for the singularity and absolute continuity of the distribution of  $\xi$ are established. It is shown that the distribution of $\xi$ is absolutely continuous if and only if $p_1=\frac{1}{3}=p_2$. In all other cases, the distribution is singular (i.e., supported on a set of Lebesgue measure zero).
When $p_0p_1p_2p_3=0$, the fractal Hausdorff–Besicovitch dimension of the spectrum (i.e., the minimal closed support) of the distribution of $\xi$ and of the essential support of its density is explicitly determined under the condition $p_{i}p_{i+1}p_{i+2}\ne0$ for $i=0,1$.

 The work also discusses the connection between the distribution of $\xi$ and infinite Bernoulli convolutions governed by the corresponding series as well as representations of numbers in the base-$3$ numeral system with one redundant digit. Several open problems are formulated.

For the numeral system with base $3$ and alphabet $A=\{0,1,2,3\}$, the problem of determining the number of representations of a number is completely solved. It is proven that almost all numbers (with respect to the Lebesgue measure) in the interval $[0;\frac{3}{2}]$ have a continuum of distinct representations, while those with a unique representation form a fractal set of Hausdorff–Besicovitch dimension $\log_3{2}$.
\end{abstract}

\maketitle

\section{Introduction}
Let $(\Omega,\mathfrak{B})$ be a measurable space and let $\mu$ and $\nu$ be two continuous probability measures defined on it.  Recall that the measure $\mu$ is called \emph{absolutely continuous} with respect to the measure $\nu$ ($\mu \leq \nu$) if $\mu(E)=0$ for all $E\in\mathfrak{B}$ such that $\nu(E)=0$.
If $\mu\leq\nu$ and $\nu \leq \mu$, the measures are called \emph{equivalent} ($\mu \sim \nu$).
The measures $\mu$ and $\nu$ are said to be \emph{orthogonal} ($\mu \perp \nu$) if there exists a set $S \in \mathfrak{B}$ such that $\mu(S) = 0$ and $\nu(\Omega \setminus S) = 0$.

If $\mu \perp \nu$ and one of the measures is the Lebesgue measure, the other is called \textit{singular}, i.e., orthogonal to the Lebesgue measure. The probability distribution function of a singularly continuous random variable has a derivative equal to zero almost everywhere with respect to the Lebesgue measure.

Singular probability measures and singular functions remain insufficiently studied objects, and their general theory is underdeveloped. Nevertheless, Zamfirescu’s theorem~\cite{Z1981} states that, in the class of continuous probability distribution functions, those that have a nontrivial singular component form a set of the second Baire category.

There are several problems related to singular measures (random variables).

1. \textit{The problem of deepening the Jessen–Wintner theorem}~\cite{jw}, which asserts that the sum of a convergent (with probability 1) random series of independent discretely distributed random variables has a pure distribution (purely discrete, purely absolutely continuous, or purely singular). The criterion for discreteness is known, but the criteria for singularity and absolute continuity are still unknown. In certain classes of random variables, this problem has been successfully solved~\cite{AGPT_2006, AGPT_2007}.

2. \textit{The problem of convolution of two singular measures (distributions)}.
Currently, the necessary and sufficient (as well as reasonably general sufficient) conditions under which the convolution of two singular measures (distributions) is singular (or absolutely continuous) are unknown.

3. \textit{The problem of convolution of two Cantor distributions} (distributions supported on nowhere dense sets of Lebesgue measure zero). To date, necessary and sufficient conditions for the spectrum (the minimal closed support) of the convolution of two Cantor distributions to be a Cantor set or Cantorval remain unknown.

4. \textit{The problem of the behavior of the absolute value of the characteristic function at infinity}.

Recall that the characteristic function of the distribution of a random variable $\zeta$ with probability distribution function $F_{\zeta}(x)$ is the complex-valued function $f_{\zeta}(t)$, which is the expected value of the random variable $e^{it\zeta}$, i.e.,
$f_{\zeta}(t)=Me^{it\zeta}=\int\limits_{-\infty}^{+\infty}e^{itx}dF_{\zeta}(x).$ It is known~\cite{p2} that for an absolutely continuous distribution of a random variable $\zeta$ with characteristic function $f_{\zeta}(t)$, the value $L_{\zeta} \equiv \limsup\limits_{|t| \to \infty} |f_{\zeta}(t)|$ is equal to $0$, for a discrete one $L_{\zeta} = 1$, and for a singular distribution $L_{\zeta} \in [0;1]$.
The value of $L_{\zeta}$ for a singular distribution characterizes its closeness to being absolutely continuous or discrete.
The nature of $L_{\zeta} = 0$ or $L_{\zeta} = 1$ for a singular distribution $\zeta$ remains mysterious.

For a continuous random variable $\zeta$ satisfying the Jessen–Wintner theorem, the condition $L_{\zeta} = 0$ implies singularity. This is the essence of the method of characteristic functions for establishing the singularity of a distribution.

  \section{Object of study}
Let $s$ and $r$ be fixed parameters such that $2\leq s\leq r$, let $A=\{0,1,...,r\}$ be an alphabet, and let $L=A\times A\times ...$ be the set of sequences of elements from the alphabet $A$. The representation of a number $x\in [0;\frac{r}{s-1}]$ by the series $x=\alpha_1 s^{-1}+\alpha_2 s^{-2}+...+\alpha_k s^{-k}+...\equiv\Delta^{r_s}_{\alpha_1\alpha_2...\alpha_k...}$ is called its representation in the base-$s$ system with a redundant alphabet. The redundancy of the alphabet generally determines the non-uniqueness of a number’s representation. Each pair of parameters  $(s;r)$ generates its own unique geometry.

Given a predetermined sequence $(\xi_k)$ of independent random variables, each taking values $0,1,...,r$ with probabilities $p_0,p_1,...,p_r$, respectively (with all $ p_i<1$), consider the distribution of the random variable $\xi=\Delta^{r_s}_{\xi_1\xi_2...\xi_k...}.$

According to the Jessen–Wintner theorem, the distribution of $\xi$ is either purely absolutely continuous or purely singular. The problem of distinguishing between these cases by necessary and sufficient conditions for arbitrary pair of parameters $s$ and $r$ is nontrivial as the problem of finding an explicit expression for the probability distribution function of $\xi$ is, even for specific values of the parameters.

With regard to the random variable $\xi$, in addition to the aforementioned issues, we are interested in the following problems:

1) In the case of singularity, to describe the topological, metric, and fractal properties of the Lebesgue zero measure set, where the distribution is supported.

2) In the case of absolute continuity, to find conditions under which $\xi$ can be decomposed into the sum of two independent components, one of which has a uniform distribution on an interval.

3) To study the behavior of the absolute value of the characteristic function of the distribution at infinity.

4) To establish the connection between the distribution of $\xi$ and infinite Bernoulli convolutions as well as convolutions of Cantor-type distributions.

This work is devoted to a partial solution of these problems.

In the paper~\cite{p8}, the case $s=2=r$ was thoroughly investigated. It was established that $\xi$ has an absolutely continuous distribution only if
$p_1= 0{,}5$ or $p_0=0{,}5=p_2$. In all other cases, the distribution of $\xi$ is singular.

In this paper, we consider the case $s=3=r$ and obtain a complete solution to the problem of determining the type of distribution of $\xi$ as well as describe the fractal properties of the distribution support under certain conditions.

In the sequel, we write
$\Delta^{r_s}_{\alpha_1\alpha_2...\alpha_k...}=
\Delta_{\alpha_1\alpha_2...\alpha_k...}=\sum\limits_{n=1}^{\infty}3^{-n}\alpha_n$, $\alpha_n\in A\equiv\{0,1,2,3\}$.
The expression $\Delta_{\alpha_1\alpha_2...\alpha_k...}=x$ is called the \emph{$\Delta$-representation} of the number $x$, and $\alpha_n$ is the $n$th digit of this representation. Parentheses in the $\Delta$-representation of a number indicate a repeating (periodic) part. A one-digit period is called simple.

\section{The connection between the distribution of $\xi$ and the infinite Bernoulli convolution governed by a series}

Every number
   $x=\Delta_{\alpha_1\alpha_2...\alpha_n...}$, $\alpha_n\in A$  is a value of the random variable $\xi$. The set $E_{\xi}$ of values of $\xi$ coincides with the set of incomplete sums (subsums) of the series
\begin{equation}\label{eq:sum1}
\sum\limits_{n=1}^{\infty}u_n=
\frac{1}{3}+\frac{1}{3}+\frac{1}{3}+
\frac{1}{3^2}+\frac{1}{3^2}+\frac{1}{3^2}+...+
\frac{1}{3^n}+\frac{1}{3^n}+\frac{1}{3^n}+...=
\sum\limits_{k=1}^{n}u_k+r_n=\frac{3}{2},
\end{equation}
that is, the set $E(u_n)=\{x: x=\sum\limits_{n=1}^{\infty}\varepsilon_n u_n,\;\varepsilon_n \in \{0;1\},\; M \in 2^N\}.$
 
Since for the terms $u_n$ and corresponding remainders $r_n$ of the given series the inequality $u_n<r_n$ holds for all $n\in N$, it follows from Kakeya's theorem~\cite{kak} that $E= [0; \frac{3}{2}]$.

The infinite Bernoulli convolution governed by the series~\eqref{eq:sum1} is the distribution of the sum
$\eta=\sum\limits_{n=1}^{\infty}u_n\eta_n,$
where $(\eta_n)$ is a sequence of independent random variables taking values $0$ and $1$ with probabilities $q_0>0$ and $q_1=1-q_0>0$, respectively.

According to the Jessen–Wintner theorem, the distribution of $\eta$ is either purely singular or purely absolutely continuous. Below, it will be shown that it is always singular.

The distribution of  $\eta$ is a special case of the distribution of $\xi$, where $p_0=q_0^3$, $p_1=3q_0^2q_1$, $p_2=3q_0q_1^2$, $p_3=q_1^3$.

If one attempts to implement a geometric approach to analyzing the distribution of $\xi$, difficulties arise due to the non-uniqueness of the $\Delta$-representation of numbers.
 
\section{The number of representations of a number and the overlapping of cylinders}

It is obvious that the equality
$\Delta_{\alpha_1...\alpha_k...}=\Delta_{c_1...c_k...}$ is equivalent to
$\sum\limits_{k=1}^{\infty}3^{-k}(\alpha_k-c_k)=0.$ 

The numbers $0=\Delta_{(0)}$ and $\frac{3}{2}=\Delta_{(3)}$ have a unique $\Delta$-representation, but they are not the only ones. There exist numbers that have a finite, countable, or even continuum set of distinct $\Delta$-representations.

Note that there are three pairs of interchangeable consecutive digit pairs in a $\Delta$-representation of a number, which do not change its value: $\overline{03}\leftrightarrow\overline{10}$,
  $\overline{13}\leftrightarrow\overline{20}$, $\overline{23}\leftrightarrow\overline{30}$.

\begin{theorem}
A number whose $\Delta$-representation contains a simple period (a period of one digit) has a countable set of different representations (except for the numbers $0$ and $\frac{3}{2}$).
 \end{theorem}
\begin{proof}
1. First, we show that the number $x_0=\Delta_{\alpha_1...\alpha_k(0)}$ has a countable set of different representations. Indeed, for an arbitrary-length sequence of twos, we have:
          $$\Delta_{\alpha_1...\alpha_k(0)}=\Delta_{\alpha_1...\alpha_{k-1}[\alpha_k-1](2)}
    =\Delta_{\alpha_1...\alpha_{k-1}[\alpha_k-1]{2...2}3(0)}.$$
   Thus, $x_0$ has an infinite set of different representations.
Let us show that this set is countable. Consider an arbitrary representation $\Delta_{\alpha_1...\alpha_k...}$
  of the number $1$.
    If we suppose $\alpha_1=0$, we get a contradiction: $1=\Delta_{0\alpha_2\alpha_3...}\leq\Delta_{0(3)}=\frac{1}{2}$.

The assumption $\alpha_1=1$ also leads to a contradiction: $\frac{1}{3}\leq 1=\Delta_{1\alpha_2\alpha_3...}\leq\Delta_{1(3)}=\frac{5}{6}$.

If $\alpha_1=3$, then $1=\Delta_{3\alpha_2...\alpha_n...}\leq \Delta_{3(0)}=1$, and thus $\alpha_k=0$ for all $k>1$, so $1=\Delta_{3(0)}=\Delta_{(2)}$.

If  $\alpha_1=2$, we analyze $\alpha_2$. 

As with $\alpha_1$, we find that $\alpha_2\notin \{0,1\}$, 
and if $\alpha_2=3$, we uniquely get $1=\Delta_{23(0)}$.

  If $\alpha_2=2$, we continue this process for $\alpha_3$, and so on. This yields
  $1=\Delta_{2...23(0)}=\Delta_{(2)}$ for all $k\in N$. Therefore, the number $1$ has a countable set of different representations.

The countability of representations of the number $\Delta_{0...0(2)}$ is justified analogously. Taking into account the equalities
     \[ \Delta_{\alpha_1...\alpha_k(0)}=\Delta_{\alpha_1...\alpha_{k-1}[\alpha_k-1](2)}=
    \Delta_{\alpha_1...\alpha_{k-1}[\alpha_k-1](0)}+\Delta_{\underbrace{0...0}_{k}(2)},\]
    \[\Delta_{\alpha_1...\alpha_k(0)}=
    \alpha_1\Delta_{1(0)}+\alpha_2\Delta_{01(0)}+...+\alpha_k\Delta_{\underbrace{0...0}_{k-1}1(0)}\]
  the proof can be considered complete.
  
2. The countability of the representations of the number  $x_0=\Delta_{c_1...c_k(2)}$ follows from item 1, since
   $\Delta_{c_1...c_k(2)}=\Delta_{c_1...c_k23(0)}$.
   
3. By replacing pairs of consecutive digits in the $\Delta$-representation (e.g., 03 with 10), we obtain:
   \[\Delta_{(1)}=\Delta_{0(3)}=\Delta_{10(3)}=\Delta_{1...10(3)}.\]
   Let $\Delta_{\alpha_1...\alpha_n...}$ be an arbitrary representation of the number $x_0=\Delta_{(1)}$.
  
Since
   $\Delta_{0(3)}=\Delta_{(1)}=x_0<\Delta_{2(0)}$, we must have $2\neq\alpha_1\ne3$. If $\alpha_1=0$, then
   $x_0=\Delta_{0(3)}$. If $\alpha_1=1$, we analyze $\alpha_2$, and so on.
   For any $k\in N$, we have
   $$\Delta_{\underbrace{1...1}_k0(3)}=\Delta_{(1)}\leq x_0< \Delta_{\underbrace{1...1}_k2(0)},$$
    so $2\neq \alpha_{k+1}\neq 3$. If $\alpha_{k+1}\neq 1$, then
    $x_0=\Delta_{\underbrace{1...1}_k0(3)}$.  
  
  A number $x_0=\Delta_{c_1...c_k(1)}$, where $c_k\neq 1$, has a countable set of different $\Delta$-representations, since the pair $c_{k}1$ has no alternative substitution, and the set $c_1...c_k$k may have only a finite
number of alternatives, or none at all if it contains no digits $0$ or $3$.   
     
  4. Since
   $\Delta_{a(3)}=\Delta_{[a+1]0(3)}=\Delta_{[a+1]1...1(3)}
   =\Delta_{[a+1](1)}$, then according to item 3, the number $x_0=\Delta_{a(3)}$ has a countable set of $\Delta$-representations. The countability of the set of representations of the number $x_0=\Delta_{c_1...c_k(3)}$ follows from item 3 and the remark above.
  Therefore, the $\Delta$-representations of the number $x_0$ are exhausted by the cases described above.   
\end{proof}

\begin{corollary}
  Almost all numbers (with respect to the Lebesgue measure) in the interval $[0;\frac{3}{2}]$ have a continuum of different $\Delta$-representations.
\end{corollary}
\begin{theorem}
  A number whose $\Delta$-representation contains only the digits $1$ and $2$, each appearing infinitely many times, has a unique $\Delta$-representation.
\end{theorem}
\begin{proof}
Let $\Delta_{\alpha_1...\alpha_n...}$ be an arbitrary  $\Delta$-representation of the number
  $$x_0=\Delta_{\underbrace{11...1}_{a_1}\underbrace{22...2}_{b_1}
  \underbrace{11...1}_{a_2}\underbrace{22...2}_{b_2}...
  \underbrace{11...1}_{a_n}\underbrace{22...2}_{b_n}...}, a_n, b_n \in N.$$

  Denote the digits of the given $\Delta$-representation of $x_0$ by $c_1$,$c_2$,...,$c_n$,...

  We shall prove that $\alpha_n=c_n$ for all natural $n$ using the principle of mathematical induction.

Since
  $\Delta_{(1)}=\Delta_{0(3)}<x_0<\Delta_{2(0)}=\Delta_{1(2)},$
it follows that $\alpha_1=1$.

Assume that $\alpha_j=c_j$ for all $j\leq k$.
Consider the digit $\alpha_{k+1}$. There are two possible cases: 1) $c_{k+1}=1$; 2) $c_{k+1}=2$. In the first case:
\[\Delta_{c_1...c_k(1)}=\Delta_{c_1...c_k0(3)}<x_0<\Delta_{c_1...c_k2(0)}=\Delta_{c_1...c_k1(2)},\]
so $\alpha_{k+1}=1=c_{k+1}$. In the second case:
\[\Delta_{c_1...c_k(1)}=\Delta_{c_1...c_k1(3)}<x_0<\Delta_{c_1...c_k3(0)}=\Delta_{c_1...c_k(2)},\]
so $\alpha_{k+1}=2=c_{k+1}$. Thus, by the principle of mathematical induction $\alpha_{n}=c_{n}$ for all $n\in {N}$. The proof for the case where the representation begins with a sequence of twos is analogous. The theorem is proven.
\end{proof}
\begin{theorem}
   If some $\Delta$-representation of a number $x$ does not contain a simple period and contains infinitely many digits $0$ and $3$, then $x$ has a continuum of different representations.
\end{theorem}
\begin{proof}
Under the assumptions of the theorem, the $\Delta$-representation of the number $x$ contains infinitely many pairs of the form $a3$ or $0b$, where
  $a\ne 3$, $b\neq 0$. Each of these pairs admits an alternative substitution. Therefore, any number satisfying the conditions of the theorem has a continuum of different $\Delta$-representations.
  \end{proof}
  \begin{corollary}
  The numbers $$x=\Delta_{\alpha_1...\alpha_k\underbrace{1...1}_{a_1}
  \underbrace{2...2}_{b_1}...\underbrace{1...1}_{a_n}\underbrace{2...2}_{b_n}...} \mbox{ and }
  y=\Delta_{\alpha_1...\alpha_k\underbrace{2...2}_{b_1}
  \underbrace{1...1}_{a_1}...\underbrace{2...2}_{b_n}\underbrace{1...1}_{a_n}...}$$ have a finite number of different $\Delta$-representations. Moreover, if  $\alpha_j=3$ for all $j=\overline{1,k}$, then both $x$ and $y$ have a unique $\Delta$-representation.
\end{corollary}
\begin{corollary}
The set of numbers that have a unique $\Delta$-representation is of continuum cardinality and has fractal Hausdorff–Besicovitch dimension equal to $\log_32$.
\end{corollary}
The set $\Delta_{c_1...c_m}$ of numbers $x\in E$ that have a $\Delta$-representation of the form $\Delta_{c_1...c_m\alpha_1...\alpha_n...}$, $(\alpha_n)\in L$, is called \textit{cylinder of rank $m$ with base $c_1...c_m$}. It is easy to prove that the cylinder  $\Delta_{c_1...c_m}$ is a segment with endpoints $a=\sum\limits_{k=1}^{m}c_k3^{-k}$,
$b=a+\frac{1}{3^m}$ and clearly,  $\Delta_{c_1...c_m}=\bigcup\limits_{i=0}^{3}\Delta_{c_1...c_mi}$.

A distinctive feature of $\Delta$-representation of numbers lies in the specific overlaps between cylinders:
\[\Delta_{c_1...c_mi}\cap\Delta_{c_1...c_m[i+1]}=\Delta_{c_1...c_{m-1}i3}=\Delta_{c_1...c_{m-1}[i+1]0},\; i=0,1,2,\]
i.e., the intersection of two ``adjacent'' cylinders is itself a cylinder of the next rank.

\section{Function with fractal level sets}
From the obtained results concerning the number of 
$r_s$-representations of numbers, one can perform a fractal analysis of the level sets of the function. Let us consider a model example with $s=3=r$.

Let us recall that the traditional quaternary representation of numbers $x\in [0;1]$ is called 
\[x=\Delta^4_{\alpha_1\alpha_2...\alpha_n...}=\frac{\alpha_1}{4}+\frac{\alpha_2}{4^2}+...+\frac{\alpha_n}{4^n}+...,\]
where $\alpha_n\in\{0,1,2,3\}$, in particular $\Delta^4_{\alpha_1\alpha_2...\alpha_k(0)}=\Delta^4_{\alpha_1\alpha_2...\alpha_{k-1}[\alpha_k-1](3)}$.

By agreeing to use, among the two quaternary representations of a quaternary-rational number, only the one that has the period (0), we define the function
 $f$ by
  \begin{equation*}
    f(x=\Delta^4_{\alpha_1\alpha_2...\alpha_n...})=\Delta_{\alpha_1\alpha_2...\alpha_n...}
  \end{equation*}

\begin{theorem}
  The function $f$  has level sets of different cardinalities:
  \begin{itemize}
    \item [1)] singletons and finite sets;
     \item [2)] countable and continuum sets,
  \end{itemize}
  including sets of positive fractal dimension.
\end{theorem}
\begin{proof}
  It is clear that the cardinality of the set $f^{-1}(y_0)$ at level $y_0=f(x_0)$ is equal to the cardinality of the set of different representations of the number $y_0\in [0;\frac{3}{2}]$.

  As noted above, the numbers $0$ and $\frac{3}{2}$ have unique representations: $f^{-1}(0)=\Delta^{4}_{(0)}$, $f^{-1}(\frac{3}{2})=\Delta^{4}_{(3)}=1$. According to the lemma, the number $\Delta_{(12)}=\frac{5}{8}$ also has a unique representation $f^{-1}(\frac{5}{8})=\Delta^4_{(12)}=\frac{2}{5}$.

   The number $y_0=\Delta_{1010(12)}=\Delta_{0310(12)}=\Delta_{0303(12)}=\Delta_{1003(12)}$ has a finite number of formally different representations, and therefore has a finite level set: $$\{\Delta^4_{1010(12)}, \Delta^4_{0310(12)}, \Delta^4_{0303(12)}, \Delta^4_{1003(12)}\}.$$

    The level set $y_0=\Delta_{(1)}=\Delta_{\underbrace{1...1}_{n}0(3)}$, as follows from the above, is countable.

   The number $y_0=\Delta_{(10)}$ has a continuum set of different representations, and thus the level set of the function $f$ is continuum:
  $y_0^{-1}=\{x: x=\Delta^4_{a_1b_1a_2b_2...a_nb_n...}\},$
  where $(a_n,b_n)\in\{(1,0),(0,3)\}$.

  Writing the number $x=\Delta^4_{a_1b_1a_2b_2...a_nb_n...}$ in the form
  \begin{align*}
    x=&\Delta^4_{a_1b_1a_2b_2...a_nb_n...}=\frac{a_1}{4}+\frac{b_1}{4^2}+\frac{a_2}{4^3}+\frac{b_2}{4^4}+...\\
    =&\frac{4a_1+b_1}{4^2}+\frac{4a_2+b_2}{4^4}+...
    =\frac{\alpha_1}{16}+\frac{\alpha_2}{16^2}+...+\frac{\alpha_n}{16^n}+...,
  \end{align*}
  where $\alpha_n\equiv4a_n+b_n\in \{3;4\}$,
  we see that the set $y_0^{-1}$ consists of points whose representations in base 16 use only two digits, 3 and 4. It is well known that the set, is a continuum, self-similar Lebesgue measure zero set with fractal dimension $\log_{16}2=\frac{1}{4}$.
\end{proof}
\section{Singularity of the distribution}

The set $E_{\xi}$ under the condition $p_0p_1p_2p_3 \ne 0$ is the spectrum of the distribution of $\xi$. If $p_0p_1p_2p_3 = 0$, then the spectrum of $\xi$ is a subset of $E$.
If either $p_0=0$ or $p_3=0$, then the spectrum is $[0;1]$ or $[\frac{1}{2};\frac{3}{2}]$, respectively. In all other cases, the spectrum is a nowhere dense fractal set. A nontrivial case arises when $p_0p_ip_3 \ne 0$ and $p_j=0$, for $i\ne j$, $j \in \{1,2\}$.
\begin{lemma}
 The set
$C[\Delta;\{0,1,3\}]=\{x:~ x=\Delta_{\alpha_1\alpha_2...\alpha_n...}, \alpha_n\in \{0,1,3\}\}$
  is a nowhere dense, $N$-self-similar set of zero Lebesgue measure, whose fractal (Hausdorff) dimension is $x=\log_{3}{\frac{3+\sqrt{5}}{2}}$.
\end{lemma}
\begin{proof}
Let $\Delta'_{c_1...c_m}\equiv\Delta_{c_1...c_m}\cap C$. Since
  $C\subset [\Delta_{3}\cup\{\Delta_{1}\cup\Delta_{0}\}], \varnothing=\Delta_{1}\cap \Delta_{3},$
  it follows that $C\stackrel{\frac{1}{3}}{\sim}\Delta'_{3}$.
 Since $\Delta'_{13}\subset C$ and $\Delta'_{13}\cap \Delta_0=\varnothing $,
  we have $C\stackrel{3^{-2}}{\sim}\Delta'_{13}$.

  Moreover, $\Delta_{00}\cup\Delta_{01}$ is congruent (isometric) to
  $\Delta_{10}\cup\Delta_{11}$, $\Delta'_{0}\cup\Delta'_{1}\stackrel{\frac{1}{3}}{\sim}\Delta'_{i0}\cup\Delta'_{i1}$,
  and \[C\subset \Delta_{3}\cup\Delta_{13}\cup[(\Delta_{00}\cup\Delta_{01})\cup
  (\Delta_{10}\cup\Delta_{11})].\]
  Furthermore, $C\stackrel{\frac{1}{3^3}}{\sim}\Delta'_{i13}$, $i=0,1$.
  Taking into account that $\Delta_{003}=\Delta_{010}$, $\Delta_{030}=\Delta_{100}$, we have
  \[\Delta'_{000}\cup\Delta'_{001}\cong \Delta'_{010}\cup\Delta'_{011}\cong\Delta'_{100}\cup\Delta'_{101}
  \cong\Delta'_{110}\cup\Delta'_{111}\stackrel{3^2}{\sim}\Delta'_0\cup \Delta'_{1},\]
  and we obtain:
  $C\subset\Delta'_3\cup\Delta'_{13}\cup[\bigcup\limits_{i_1=0}^{1}...\bigcup\limits_{i_n=0}^{1}\Delta'_{i_1...i_n13}]\cup D,$
   where $D$ is a countable set, and $C\stackrel{k}{\sim}\Delta'_{i_1...i_n13}$, with $k=3^{-(n+2)}$.
   Hence, $C$ is an $N$-self-similar set, and its dimension is the solution of the equation
  $\frac{1}{3^x}+\sum\limits_{n=0}^{\infty}\frac{2^n}{3^{(n+2)x}}=1$, that is
  $x=\log_{3}\frac{3+\sqrt{5}}{2}$.
\end{proof}
\begin{remark}
The set $C[\Delta; \{0,2,3\}]$  has the same properties as $C[\Delta;\{0,1,3\}]$. The proof is analogous to the one above.
\end{remark}
\begin{theorem}
  Let $p_0p_1p_2p_3=0$.  1. If $p_i=p_{i+1}=p_{i+2}=\frac{1}{3}$, then the random variable $\xi$ has a \emph{uniform distribution} on the unit interval.

  2. If there exists a $p_j$ such that $0\neq p_j\neq \frac{1}{3}$, then the distribution of the random variable $\xi$ is a \emph{singular distribution}, and:

 2.1) if two of the probabilities are zero, then $\xi$ has a \emph{Cantor-type distribution} with a spectrum of fractal dimension $\log_{3}2$;

 2.2) if only one of the probabilities is zero, and $p_1p_2=0$, then $\xi$ has a \emph{Cantor-type distribution} with a spectrum of fractal dimension $\log_{3}\frac{3+\sqrt{5}}{2}$;

  2.3) if only one of the probabilities is zero, and $p_0p_3=0$, then $\xi$ has a \emph{singular distribution with a strictly increasing probability distribution function}, and the fractal dimension of the essential support of its density is $-\log_{3}{p_1^{p_1}p_2^{p_2}p_i^{p_i}}$.
\end{theorem}
\begin{proof}
1. If $p_3=0$, then $\xi$ is a random variable whose digits ${\xi_n}$ in the ternary representation take values  $0,1,2$ with probability $\frac{1}{3}$. Such a random variable is known to have a uniform distribution on the interval $[0;1]$.

If $p_0=0$, then consider the random variable $\hat{\xi}=\xi-\frac{1}{2}=\xi-\sum\limits_{n=1}^{\infty}\frac{1}{3^n}$, i.e., the random variable $\sum\limits_{n=1}^{\infty}\frac{\xi_n-1}{3^n}=\Delta^3_{\xi_1'\xi_2'...\xi_n'...}$, where $\xi_n'$ takes values $0$, $1$, $2$  with equal probabilities $\frac{1}{3}$. Since $\xi$ and $\hat{\xi}$ have equivalent distributions, and $\hat{\xi}$ is uniformly distributed, it follows that $\xi$ is uniformly distributed on $[\frac{1}{2};\frac{3}{2}]$.

2.1) If $p_3=0=p_i$, then the spectrum of the distribution of $\xi$ is a Cantor-type set $C[r_s, V]$, where $V=\{0,1,2\}\setminus \{i\}$, which is known to be a self-similar set with fractal dimension $\log_32$. A similar situation occurs when $p_1=0=p_2$.
 If $p_3\ne 0= p_0$ and $p_1p_2=0,$ then the corresponding random variables $\xi$ and $\hat{\xi}$ have equivalent distributions. This last case corresponds to a random variable with independent digits in the classical ternary representation. Since one digit has zero probability, $\xi$ has a singular Cantor-type distribution.

2.2). In the case $p_1p_2=0$, the spectrum of the distribution of $\xi$ coincides with the set $C[r_s; \{0,2,3\}]$ or $C[r_3;\{0,1,3\}]$. Therefore, Claim 2.2) follows from the previous lemma and remark.

 2.3). Let $p_0p_3=0$. If $p_3=0$, then the random variable $\xi$ has independent and identically distributed ternary digits, and if $p_j\ne \frac{1}{3}$ for some $j$, then the distribution is singular.

Indeed, in this case, the set of growth points of the probability distribution function is the interval $[0;1]$. The probability that $\xi$  belongs to the set  $E[3; p_0,p_1,p_2]$ of numbers whose digit frequencies match $p_i$ for each $i=0,1,2$ is equal to 1. As it was shown by Eggleston~\cite{Eggleston H.G.} and Billingsley~\cite{b}, the set $E$ has Lebesgue measure zero and Hausdorff–Besicovitch dimension $-\log_{3}{p_0^{p_0}p_1^{p_1}p_2^{p_2}}$. Therefore, $\xi$ has a singular distribution with a strictly increasing distribution function~\cite{Marsaglia1983}.
 
If $p_0=0$, then the same conclusion holds for the random variable $\hat{\xi}=\xi-\frac{1}{2}$, whose distribution is equivalent to distribution of $\xi$.
\end{proof}
In the case $p_0p_1p_2p_3\neq 0$, we use the method of characteristic functions.

The characteristic function of the distribution of $\xi$
has the form
\[f_{\xi}(t)=\prod\limits_{k=1}^{\infty}\varphi_k(t),\]  where  $\varphi_k(t)=p_0+p_1e^{it\cdot3^{-k}}+p_2e^{i2t\cdot3^{-k}}+p_3e^{i3t\cdot 3^{-k}}=
 \sum\limits_{m=0}^{3}p_m\cos{\frac{mt}{3^k}}+i\sum\limits_{m=0}^{3}p_m\sin{\frac{mt}{3^k}},$
 and satisfies the functional equation $f_{\xi}(t)=\varphi_1(t)f_{\xi}(\frac{t}{3})$. In particular, $f_{\xi}(2\pi ns)=f_{\xi}(2\pi n)$, $\forall n\in N$ and
 $L_{\xi}\geq|f_{\xi}(2\pi)|$.

\begin{theorem}\label{th:s}
  If $p_1=\frac{1}{3}=p_2$, then the random variable $\xi=\Delta_{\xi_1\xi_2...\xi_n...}$ with independent digits $\xi_n$ taking values $0,1,2,3$ with probabilities $p_0$, $p_1$, $p_2$, $p_3$, respectively, has an absolutely continuous distribution. In this case, the distribution of $\xi$
  is the convolution of a uniform distribution on $[0;1]$ and a singular Cantor-type distribution. In all other cases, the distribution of $\xi$ is singular.
\end{theorem}
\begin{proof}
If $p_1=\frac{1}{3}=p_2$, then $\xi$ can be expressed as the sum of two independent random variables: $\tau$, which has a uniform distribution on $[0;1]$, and
 $\zeta$, which has a singular Cantor-type distribution.

Indeed, a uniformly distributed random variable on $[0;1]$ has the form $\tau=\sum\limits_{n=1}^{\infty}3^{-n}\tau_n$, where
  $(\tau_n)$ are i.i.d. random variables taking values $0$, $1$, $2$ with equal probabilities $\frac{1}{3}$. The random variable $\zeta=\sum\limits_{n=1}^{\infty}3^{-n}\zeta_n$, where $\zeta_n$ are independent random variables taking values $0$, $1$ with probabilities $x>0$ and $1-x>0$, has a singular Cantor-type distribution.

The convolution of the distributions of $\tau$ and $\zeta$ is the distribution of sum $\tau+\zeta$ if addends are independent. The identity $\xi=\tau+\zeta$ is equivalent to system of equations $\xi_n=\tau_n+\zeta_n$, $n=1,2,...$. This system is equivalent to the following consistent system of equations: 
  \begin{align*}
    p_{0}=P\{\xi_n=0\}=&\frac{1}{3}x,\;\;
     p_1=P\{\xi_n=1\}=\frac{1}{3}(1-x)+\frac{1}{3}x=\frac{1}{3},\\
     p_2=P\{\xi_n=2\}=&\frac{1}{3}x+\frac{1}{3}(1-x)=\frac{1}{3},\;
     p_3=P\{\xi_n=3\}=\frac{1}{3}(1-x),
  \end{align*} which implies $p_1=\frac{1}{3}=p_2$, $p_0+p_3=\frac{1}{3}$, and $x=3p_0$.
 
Thus, under the condition $p_1=\frac{1}{3}=p_2$ the distribution of $\xi$ is a convolution of uniform and singular distributions, hence it is absolutely continuous.

Now consider the case where $p_1\ne \frac{1}{3}$ or $p_2\ne \frac{1}{3}$. We use the characteristic function method to obtain sufficient conditions for singularity, specifically proving that $L_{\xi}>0$, which implies singularity. We first show that the infinite product
\begin{equation}\label{d}
  \prod\limits_{k=2}^{\infty}\varphi_k(2\pi)\ne 0.
\end{equation}

Since $\sum\limits_{m=1}^{3}p_m\sin{\frac{2 \pi m}{3^k}}>0$  for all $k>1$, it follows that $\varphi_k(2\pi)\ne 0$.
It is known that the absolute convergence of the infinite product~\eqref{d} is equivalent to the convergence of the series 
$\sum\limits_{k=2}^{\infty}|\varphi_k(2\pi)-1|$. We prove the convergence of this series using the comparison test. Let us estimate $u_k=|\varphi_k(2\pi)-1|$:
\begin{align*}\label{eq:varphi}
u_k &= \left[(\sum\limits_{m=0}^{3}p_m\cos(2\pi m3^{-k})-1)^2+
(\sum\limits_{m=0}^{3}p_m\sin(2\pi m3^{-k}))^2\right]^{\frac{1}{2}}=\\
  &= \left[(\sum\limits_{m=0}^{3}p_m)^2+
  4(\sum\limits_{m=1}^{3}p_m\sin^2(\pi m 3^{-k})-
  \sum\limits_{0\neq i\neq j\neq 0}p_ip_j\sin^2(\pi(i-j)3^{-k}))-1\right]^{\frac{1}{2}}=\\
  &= 2\left[\sum\limits_{m=0}^{3}p_m\sin^2(\pi m3^{-k})-\sum\limits_{0\neq i\neq j\neq 0}p_ip_j\sin^2(\pi(i-j)3^{-k})\right]^{\frac{1}{2}}\leq\\
  &\leq 2\left[\sum\limits_{m=1}^{3}p_m\sin^2(\pi m3^{-k})\right]^{\frac{1}{2}}\leq
  2\left[\sum\limits_{m=1}^{3}p_m\sin^2(3\pi\cdot3^{-k})\right]^{\frac{1}{2}}\leq
  2\sin(3\pi\cdot 3^{-k})\leq 6\pi\cdot3^{-k}.
\end{align*}

Let us note that the last three inequalities are valid for all  $k>2$. From the final inequality and the comparison test, it follows that the series $\sum\limits_{k=2}^{\infty}|\varphi_k(2\pi)-1|$ converges, and hence so does the infinite product~\eqref{d}, which means that inequality~\eqref{d} holds. Therefore, $|f_{\xi}(2\pi)|=0$ if and only if $\varphi_1(2\pi)=0$, i.e., when $\sum\limits_{m=0}^{3}\cos{\frac{2\pi m}{3}}=0=\sum\limits_{m=0}^{3}\sin{\frac{2\pi m}{3}}$, which is equivalent to the system of equations
$p_0-\frac{1}{2}p_1-\frac{1}{2}p_2+p_3=0=\frac{\sqrt{3}}{2}(p_1-p_2).$ Thus, $|f_{\xi}(2\pi)|=0$ if and only if $p_1=\frac{1}{3}=p_2$, and under this condition, as shown earlier, the distribution of $\xi$ is absolutely continuous.

Hence, when $p_1\neq\frac{1}{3}$ or $p_2\neq \frac{1}{3}$, we have $|f_{\xi}(2\pi)|\ne 0$, and therefore $L_{\xi}>0$ which implies that the distribution of $\xi$ is singular.
\end{proof}
\begin{corollary}
The random variable $\eta$ has a purely singular distribution.
\end{corollary}
  The distribution of $\eta$ is singular because the system $3q_0^2(1-q_0)=\frac{1}{3}=3q_0(1-q_0)^2$ is inconsistent.
\section{The distribution of $\xi$  as a convolution of Cantor-type distributions}
\begin{lemma}
The sum of two independent singularly distributed random variables   $\theta=\sum\limits_{n=1}^{\infty}\theta_n3^{-n}$ and $\varepsilon=\sum\limits_{n=1}^{\infty}\varepsilon_n3^{-n},$
where $(\theta_n)$ and $(\varepsilon_n)$ are sequences of independent random variables taking values $0$, $2$ and $0$, $1$, respectively, with probabilities $u$, $1-u$ and  $v$, $1-v$ has a singular distribution whose spectrum is the interval $[0;\frac{3}{2}]$.
\end{lemma}
\begin{proof}
  Indeed, the distribution of the sum $\theta+\varepsilon$ belongs to the same class of distributions as $\xi$ with corresponding parameters $p_0=uv$,
$p_1=u(1-v)$, $p_2=v(1-u)$, $p_3=(1-u)(1-v)$. However, the system $u(1-v)=\frac{1}{3}=v(1-u)$ is inconsistent. Therefore, by Theorem~\ref{th:s}, the distribution of $\theta+\varepsilon$ is singular.
\end{proof}

\begin{theorem}
  If $p_0=(p_0+p_1)(p_0+p_2)$, then the distribution of $\xi$ is the convolution of two Cantor-type distributions, namely, the distributions of the random variables 
   $\theta=\Delta_{\theta_1\theta_2...\theta_n...}, \varepsilon=\Delta_{\varepsilon_1\varepsilon_2...\varepsilon_n...},$
where $(\theta_n)$ and $(\varepsilon_n)$ are sequences of independent random variables taking values $0$, $2$ and $0$, $1$ with probabilities $p_0+p_1$, $1-(p_0+p_1)$ and  $p_0+p_2$, $1-(p_0+p_2)$, respectively.
\end{theorem}
\begin{proof}
   If $P\{\theta_n=0\}=x$, $P\{\theta_n=2\}=1-x$, $P\{\varepsilon_n=0\}=y$, $P\{\varepsilon_n=1\}=1-y$, then the equality
  $\xi=\theta+\varepsilon$ is equivalent to the system $p_0=xy$, $p_1=x(1-y)$, $p_2=(1-x)y$, $p_3=(1-x)(1-y)$,
  which, in turn, is equivalent to $x=p_0+p_1$, $y=p_0+p_2$, $p_0=xy=(p_0+p_1)(p_0+p_2)$.
\end{proof}
\begin{remark}
   For $p_0=p_1=p_2=p_3=\frac{1}{4}$, the condition $p_0=(p_0+p_1)(p_0+p_2)$ is satisfied. In this case, the distribution of $\xi$ is singular.
\end{remark}
\section{Concluding remarks}
 We were not able to obtain an explicit expression for the probability distribution function $F_\xi(x)$. The problem of describing the topological, metric, and fractal properties of the essential support of the density of distribution of $\xi$, $N_\xi = \{x: F'_{\xi}(x)>0 \mbox{ or } F'_{\xi}(x) \mbox{ does not exist}\}$, under the condition $p_0p_1p_2p_3\ne 0$ remains open.
The main difficulties arise from the non-uniqueness of the number representation by the series and the overlapping of cylinders. Despite the relatively simple structure of the overlaps, computing the probabilities of cylinders remains challenging.

\end{document}